\documentclass[12pt]{article}
\usepackage{amssymb}
\def\no{\noindent}

\def\f#1#2{\frac{#1}{#2}}
\def\f{\frac}

\def\ov{\overline}

\def\R{{\mathbb R}}
\def\S{{\cal S}}
\def\D{{\cal D}}
\def\C{{\cal C}}

\def\ep{\epsilon}

\def\p{\partial}

\def\L{{\bf L}}

\def\N{{\cal N}}

\def\v{\vskip 1em}

\def\ve{\varepsilon}

\def\meas{\hbox{meas}\,}

\def\T{{\cal T}}
\def\tv{\hbox{Tot.Var.}}
\def\Tilde{\widetilde}

\def\sqr#1#2{\vbox{\hrule height .#2pt
\hbox{\vrule width .#2pt height #1pt \kern #1pt
\vrule width .#2pt}\hrule height .#2pt }}
\def\square{\sqr74}
\def\endproof{\hphantom{MM}\hfill\llap{$\square$}\goodbreak}

\newcommand{\beq}{\begin{equation}}
\newcommand{\eeq}{\end{equation}}
\newcommand{\ben}{\begin{eqnarray}}
\newcommand{\een}{\end{eqnarray}}
\newcommand{\beno}{\begin{eqnarray*}}
\newcommand{\eeno}{\end{eqnarray*}}
\newtheorem{Theorem}{Theorem}[section]

\newtheorem{Def}{Definition}[section]

\newcommand{\QED}{$\square$}
\begin{document}
\title{\bf  On Asymptotic Variational Wave Equations}
\author{Alberto Bressan$ ^1$, Ping Zhang$ ^2$, and Yuxi Zheng$ ^1$ \\[2mm]
{\small $ ^1$ Department of Mathematics, Penn State University,  
PA 16802.}\\[-2mm]
{\small E-mail: bressan@math.psu.edu; yzheng@math.psu.edu}\\[-2mm]
{\small $ ^2$ Academy of Mathematics and System Sciences, CAS,
Beijing 100080,}\\[-2mm]
{\small  China. \qquad  E-mail: zp@mail.math.ac.cn} }
\maketitle

\begin{abstract}
We investigate the equation $(u_t + (f(u))_x)_x = f''(u) (u_x)^2/2$
where $f(u)$ is a given smooth function. Typically $f(u)= u^2/2$ or $u^3/3$.
This equation models unidirectional and weakly nonlinear waves for
the variational wave equation
$u_{tt} - c(u) (c(u)u_x)_x =0$ which models some liquid crystals
with a natural sinusoidal $c$. The equation itself is also the 
Euler-Lagrange
equation of a variational problem.
Two natural classes of solutions can be associated with this equation.
A conservative solution will preserve its energy in time, while a
dissipative weak solution loses energy at the time when
singularities appear.
Conservative solutions are globally defined,
forward and backward in time, and preserve interesting geometric features,
such as the Hamiltonian structure. 
On the other hand, dissipative solutions appear to be more natural
from the physical point of view.

We establish the well-posedness of the Cauchy problem
within the class of conservative solutions,
for initial data having finite energy and assuming that the 
flux function $f$ has
Lipschitz continuous
second-order derivative. In the case where $f$ is convex,
the Cauchy problem is well-posed also within the class of
dissipative solutions.
However, when $f$ is not convex,
we show that the dissipative solutions
do not depend continuously on the initial data.

\end{abstract}
\vspace*{3mm}
\renewcommand{\theequation}{\thesection.\arabic{equation}}
\setcounter{equation}{0}
\section{Introduction}

A nonlinear variational wave equation whose
wave speed is a sinusoidal function of the wave amplitude arises in the
study of nematic liquid crystals. It is given by
\beq
\p_{t}^2\psi-c(\psi)\p_x(c(\psi)\p_x\psi)=0,
\label{1.1}
\eeq
with
\beq
c^2(\psi) = \alpha\sin^2(\psi) + \beta\cos^2(\psi),
\label{1.5}
\eeq
where $\alpha$ and $\beta$ are positive physical constants.
We refer the reader to \cite{[14]}, \cite{[15]}, \cite{[20]} for
background information on the equation.
Glassey, Hunter, and Zheng \cite{[13]} have shown that singularities
can form from smooth data for equation (\ref{1.1})-(\ref{1.5}).
Assuming that the wave speed  $c(\cdot)$ 
is a monotone increasing function,
the global existence of (dissipative)
weak solutions has been established in
\cite{[26]}, \cite{[27]}, \cite{[28]}, \cite{[29]}.
The general problem of the global existence and uniqueness of
conservative solutions to the Cauchy problem of equation (\ref{1.1}) 
will be addressed in a forthcoming paper \cite{[BZ]}.

The study of solutions to
(\ref{1.1})-(\ref{1.5}) 
consisting of a small-amplitude and high-frequency perturbation of 
a constant state has greatly contributed to
the understanding of this equation \cite{[13]}, \cite{[26]},
\cite{[27]}, \cite{[28]}, \cite{[29]}.
Hunter and Saxton first studied these waves in \cite{[15]}.
Given a constant state $a$, these perturbed solutions
take the form
$$
\psi(t,x ) = a + \ep u(\ep t, x-c(a)t) + O(\ep^2). 
$$
Hunter and Saxton found that $u(\cdot, \cdot)$ satisfies
\beq
(u_t+u^nu_x)_x = \f12n u^{n-1}u^2_x \label{1.7}
\eeq
up to a scaling and reflection of the independent variables,
assuming that $a$ is such that $c^{(k)}(a) =0, k=1, 2, \dots n-1$,
but $c^{(n)}(a) \ne 0$, for an integer $n\ge 1$.
In connection with 
our sinusoidal function $c$ modeling 
nematic liquid crystals in (\ref{1.5}), the relevant approximations in
(\ref{1.7}) are those with $n=1,2$. 
The case $n=1$ yields the {\it first-order asymptotic equation}
\beq
(u_t + uu_x)_x = \f12u_x^{\,2},  \label{1.10}
\eeq
for which existence and uniqueness of admissible conservative and
dissipative weak solutions have both
been established, see \cite{[16]} and \cite{[23]}, \cite{[24]},
\cite{[25]}.
This equation is also an
asymptotic equation of the Camassa-Holm equation \cite{[3]},
describing the motion of
solitary waves in shallow water. 
For recent literature on the Camassa-Holm equation,
we refer the reader to \cite{[4]}, \cite{[5]}, \cite{[6]}, 
\cite{[7]}, \cite{[8]}, \cite{[22]}
and in particular \cite{[BC]}.

The case $n=2$ yields the {\it second-order asymptotic equation}
\beq
(u_t+u^2u_x)_x = uu^2_x\,.
\label{1.15}
\eeq
In  \cite{zzpV} Zhang and Zheng established
that dissipative solutions exist for (\ref{1.15}) with BV data.
In the analysis of (\ref{1.1}), a major difficulty
is concentration of energy at points 
where $c'=0$, as in the example on p.~70
of \cite{[13]}. We hope that investigation of singularities of the same type
for the second-order
asymptotic equation will be helpful toward the understanding of the
original equation (\ref{1.1}).

Rather than (\ref{1.7}), in the present paper we study
a somewhat more general class of equations, having the form 
\beq
(u_t+f(u)_x)_x=\f12 f''(u)u_x^2\,. \label{1.16}
\eeq
Here $u=u(t,x)$ is a scalar function
defined for $(t,x)\in \R_+\times\R_+$ where $\R_+\doteq [0,\infty[\,$, and 
$f$ is a $C^2$ function. More restrictions on $f$ will
be specified later. As initial and boundary data
we take
\beq
u(0,x)=\bar u(x)\,,\qquad\qquad u(t, 0)=0\,.\label{1.17}
\eeq
Integrating equation (\ref{1.16}) w.r.t.~$x$, we obtain
\beq
u_t+f(u)_x=\f12 \int_0^x f''(u)u_x^2\,dx \,.\label{1.18}
\eeq
It is now clear that,
to make sense of this equation, we should require that the
function 
$u(t,\cdot)$ be absolutely continuous with  
derivative $u_x(t,\cdot)$ locally square integrable, for every 
fixed time $t$.
Moreover, to satisfy the boundary condition at $x=0$,
one needs the nonnegativity of the characteristic speed at $u=0$,
namely 
\beq
f'(0)\geq 0. \label{1.20}
\eeq

The local integrability assumption $u_x(t,\cdot)\in \L^2_{loc}(\R_+)$
imposes 
a certain degree of regularity on the function $u$. 
Therefore, there is no 
need to consider weak solutions in distributional sense
and a stronger concept of solution can be adopted.

\begin{Def}{\label{Definition 1.}}

\no A function $u=u(t,x))$ defined 
on $[0,T]\times\R_+$ is a {\bf solution} of the initial-boundary
value problem (\ref{1.17}) --(\ref{1.20}) if the following holds.

\smallskip
\no (i) 
The function $u$ is locally 
H\"older continuous w.r.t~both variables $t,x$.
The initial and boundary conditions (\ref{1.17}) hold pointwise.
For each time $t$, the map 
~$x\mapsto u(t,x)$ is absolutely continuous with 
$u_x(t,\cdot)\in \L^2_{loc}
(\R_+)$.

\no (ii) For any $M>0$, consider 
the restriction of $u$ to the interval $x\in [0,M]$.
Then the map $t\mapsto u(t,\cdot)\in \L^2([0,M])$
is absolutely continuous and satisfies the equation
\beq 
\frac{d}{dt} u(t,\cdot)= -f'(u) u_x+\f12 \int_0^* 
f''(u)\,u_x^2\,dx \label{1.21}
\eeq
for a.e. ~$t\in[0,T]$.  Here equality is understood 
in the sense of functions in $\L^2([0,M])$.
\end{Def}

In spite of the regularity assumptions, the requirements contained in
the above definition are still not enough to single out a unique solution. 
Let us consider a simple example.

\bigskip
\noindent{\bf Example 1.}
Consider the flux $f(u)=u^2$ and choose the initial data
$$
u(0,x)=\left\{
\begin{array}{ll}
 -x, & \qquad 0 \leq  x \leq  1, \nonumber \\
  -1, &  \qquad x > 1. \nonumber
\end{array}
\right. 
$$
For $t\in [0,1[\,$, a  solution to (\ref{1.18}) can be constructed
by the method of characteristics, namely
$$
u(t,x)=\left\{
\begin{array}{ll}
 -x/(1-t), & \qquad 0 \leq  x \leq  (1-t)^2, \nonumber \\
  -(1-t), &  \qquad x \geq (1-t)^2. \nonumber
\end{array}
\right. 
$$
Notice that the norm of the gradient $\|u_x(t)\|_{L^\infty}$
blows up as $t\to 1$.
For $t=1$ we have $u(1,x)= 0$ for all $x\geq 0$.
At this stage, there are infinitely many ways to
further prolong the solution.
For example, we could set
\beq
u(t,x)\equiv 0\qquad\qquad t\geq 1\, ~~x\geq 0\,. \label{1.22}
\eeq
Or else we could choose an arbitrary point $b\geq 0$, an arbitrary amount of
energy $k>0$ and a time $\tau\geq 1$ and define
$$u(t,x)=0\qquad \hbox{for} ~~1\leq t\leq\tau\,,$$
while, for $t>\tau$,
$$
u(t,x)=\left\{
\begin{array}{ll} 0 & \qquad 0\leq x\leq b, \nonumber \\
 (x-b)/(t-\tau) & \qquad b \leq  x \leq  k (t-\tau)^2 +b, \nonumber \\
  k(t-\tau) &  \qquad x > b+ k(t-\tau)^2. \nonumber
\end{array}
\right.
$$
Among all these solutions, two in particular can be singled out.
If we insist that the future configurations $u(t,\cdot)$ for $t>1$
should be entirely determined only by the present configuration 
$u(1,\cdot)$,
then the only reasonable choice is (\ref{1.22}).
On the other hand, if we look for solutions that
satisfy the additional conservation equation
$$
(u_x^2)_t+(2u \,u_x^2)_x=0\,,
$$
the natural choice should be
$$
u(t,x)=\left\{
\begin{array}{ll}
x/(t-1), & \qquad 0 \leq  x \leq  (t-1)^2 \nonumber \\
   t-1,  &  \qquad x \geq (t-1)^2 \nonumber
\end{array}
\right. \qquad t>1\,.
$$
To express the fact that at time $t=1$ this solution is
different from the null solution, in some way we should think
its derivative $u_x$ as being not the zero function but
the square root of a Dirac distribution at the origin. \QED

\medskip

In the following, we say that a solution $u=u(t,x)$ is 
{\bf dissipative} if the family of absolutely continuous measures
$\{\mu_{(t)}~;~~t\geq 0\}$ defined by $d\mu_{(t)}=u_x^2(t)\,dx$
provides a measure-valued solution to
\beq 
w_t+(f'(u) w)_x\leq 0.
\eeq
More precisely, we require that
\beq
\int \phi(t,\cdot) \,d\mu_{(t)}\Big|_{t_1}^{t_2}~\leq~
 \int_{t_1}^{t_2} \left[\int [\phi_t(t,\cdot)+
\phi_x(t,\cdot)\,f'(u(t,\cdot))]\,d\mu_{(t)}\right] 
\,dt
\label{1.23}
\eeq
for {\it every} $t_2>t_1\geq 0$ and any function $\phi\in \C^1_c$,
$\phi\geq 0$. 

On the other hand, to define a semigroup of conservative solutions we
need to consider a domain $\D$ of couples $(u,\mu)$ where
$u:\R_+\mapsto\R$ 
is an absolutely continuous function with square integrable derivative
and $\mu$ is a nonnegative measure on $\R_+$.
Decomposing $\mu=\mu^a+\mu^s$ as a sum of an absolutely continuous 
and a singular part (w.r.t.~Lebesgue measure), we shall require that
$d\mu^a=u_x^2\,dx$.   We say that 
a map $t\mapsto (u(t),\,\mu_{(t)})\in\D$
is a {\bf conservative solution}
of (\ref{1.17})-(\ref{1.20}) if $u$ is a solution according to Definition 1.1
and (\ref{1.23}) is satisfied as an equality, for
{\it all} $t_2>t_1\geq 0$ and $\phi\in \C^1_c$.

As mentioned earlier, Zhang and Zheng have established in \cite{zzpV} 
the finite-time singularity formation in smooth solutions 
and the global existence of a dissipative weak solution to 
(\ref{1.17})--(\ref{1.20}) with initial data ${\bar u(x)}$ 
whose derivative is in $BV$, for $f(u) = u^3/3$.

In the present paper, we consider a general flux $f$ with
Lipschitz continuous second-order derivative such that $f'(0)\geq 0$.  
The initial data are chosen in the set of absolutely
continuous functions $\bar u$, with $\bar u(0)=0$ and
$\bar u_x\in \L^2(\R_+)$.
Our
main results can be summarized as follows.

\v
\noindent {\bf 1.} A flow of conservative 
solutions can be globally defined, forward and backward in time
(Theorem 3.1).
The conservative solution of the initial-boundary value problem
(\ref{1.17})--(\ref{1.20}) is unique, provided that a suitable 
non-degeneracy condition is satisfied (Theorem 4.1).
\v
\noindent {\bf 2.} Assuming, in addition, that the flux $f$ is convex,
then there also exists a continuous semigroup of dissipative solutions.
The dissipative solution of the initial-boundary value problem
(\ref{1.17})--(\ref{1.20}) is unique (see Theorem \ref{Theorem5.1}).
\v
\noindent {\bf 3.} If the flux $f$ is not convex, the dissipative solutions
do not depend continuously on the initial data, in general
(see Example 2 in Section 6).
\v

Before proving the main results in Section 3,
we briefly discuss the action principle and
some admissibility conditions,
whose aim is to identify a unique 
physically relevant solution to equations (\ref{1.17})--(\ref{1.20}).

\bigskip

\setcounter{equation}{0}
\section{Remarks on admissibility conditions}

The decay estimate
$$
u_x(t,x)\leq 2/t
$$
was used as an admissibility  
criterion for dissipative solutions of the first-order
asymptotic equation in \cite{[16]}, \cite{[23]}, \cite{[24]}
and \cite{[25]}.  We remark, however,  that this
is not appropriate in connection with
dissipative solutions of (1.6). Indeed, for
a solution of the second-order asymptotic equation,
the gradient $u_x$ can approach $+\infty$ as well as $-\infty$.

Another common criteria for selecting physically admissible solutions
is by vanishing viscosity. 
One might conjecture that dissipative solutions are precisely the limits
of vanishing viscosity approximations. 
We believe this is indeed the case when the flux function $f$ is 
convex, see some proofs in \cite{[16]} \cite{[22]} for $f=u^2/2$.
On the other hand, when $f$ is not convex, the dissipative solutions
do not depend continuously on the initial data (see Section 6).
We observe that the set of vanishing viscosity
limits is closed, connected and depends on the initial data in an upper
semicontinuous way. Therefore, by a topological argument, 
the vanishing viscosity criterion cannot single out a unique
limit, in general.

Concerning the vanishing dispersion limit,
numerical experiments performed with a
convex $f$ seem to indicate that
vanishing dispersion selects the conservative solutions,
see \cite{[16]}.

Next, we discuss the admissibility of solutions in terms of a
variational principle.
For all asymptotic equations (\ref{1.7}) we have the action functionals
\beq
{\cal A}_n\doteq
\int_{t_1}^{t_2}\int [u_xu_t + u^n (u_x)^2]\, dx\,dt\ . \label{2.2}
\eeq
In other words, the Euler-Lagrange equations  satisfied by
functions $u$ that render stationary the action ${\cal A}_n$
are precisely the 
asymptotic equations (\ref{1.7}). These can be derived from the
nonlinear variational wave equation (\ref{1.1})
\beq
\psi_{tt} - c(\psi)(c(\psi)\psi_x)_x =0  \label{2.3}
\eeq
by a perturbation argument.  Notice that
(\ref{2.3}) is the Euler-Lagrange equation corresponding to the Lagrangean
\beq
{\cal L} = \psi_t^2 - c^2(\psi)\psi_x^2.
\eeq
This arises often in physical models.
For weakly nonlinear waves of the form
$$
\psi = \psi_0 + \epsilon u (\tau, \theta) +
\epsilon^2 v(\tau, \theta) +O(\epsilon^3)
$$
with
$$
\tau = \epsilon t\,, \qquad \theta = x- c_0t\,, \qquad 
c_0 \doteq c(\psi_0)\,,
$$
assuming that $c_0' \doteq c'(\psi_0)\ne 0$
we have 
$$
\psi_{tt} - c(\psi)(c(\psi)\psi_x)_x = -2c_0\epsilon^2
\Big\{ (u_\tau +c_0'uu_\theta)_
\theta -\frac{1}{2}c_0'u_\theta^2\Big\} + O(\epsilon^3)\,.
$$
Moreover
$$
\psi_t^2 - c^2(\psi)\psi_x^2 = -2c_0\epsilon^3
[u_\tau u_\theta + c_0' u u_\theta^2] + O(\epsilon^4).
$$
Therefore, $u $ satisfies the first order asymptotic equation.
The corresponding Lagrangean, approximated to order $O(\epsilon^3)$, is
$-(u_\tau u_\theta + c_0' u u_\theta^2)$.

At first sight, 
one might hope that the physically relevant solutions to
the equations (\ref{1.7}) are those which maximize (or minimize)
the action in (\ref{2.2}).  Unfortunately this is not the case,
because the action ${\cal A}_n$ is not coercive.  For any
smooth solution $u$ of (\ref{1.7}) one can find compactly
supported perturbations $u+\epsilon v$ which increase the
value of ${\cal A}_n$, and others which decrease it.
The extremality of the action thus cannot be used as a selective criterion.


\setcounter{equation}{0}

\section{Conservative solutions}

We consider the evolution problem described by the equation
\beq
u_t + f(u)_x = \f12\int_0^x f''(u)u_x^2\, dx\qquad\qquad
\hbox{for all}~~t\geq 0\,,~~ x\geq 0,
\label{3.1}
\eeq
together with the boundary conditions
\beq
u(t, 0)=0\qquad\qquad \hbox{for all}~~ t\geq 0\,.
\label{3.2}
\eeq
We assume that $f \in C^2(\R)$ and
\beq
f'(0) \geq 0, \qquad 
|f''(u) - f''(v) | \leq L |u -v|, \quad \forall u, v \in \R 
\label{L0}
\eeq 
for a constant $L$.

One easily checks that every smooth solution
satisfies the additional 
conservation law for the ``energy'' $u_x^2$, namely
\beq
(u_x^2)_t+\left[ f'(u) (u_x^2)\right]_x=0\,.
\label{3.3}
\eeq
It is therefore natural to seek a continuous flow associated with
(\ref{3.1})-(\ref{3.2}) which preserves the energy 
$\int_0^\infty  u_x^2(t, x)\,dx$.
However, Example 1 in the Introduction already pointed
out a basic difficulty
which one encounters while constructing a
semigroup in the space $H^1_{\rm loc}\,$.
Indeed, when the gradient $u_x$ blows up, 
all the energy is concentrated at a single point, so that
the measure $u_x^2\,dx$ approaches a Dirac mass.

Motivated by this example, to the equations (\ref{3.1})-(\ref{3.2})
we will associate an evolution semigroup on a domain
$\D$ defined as follows.
An element of $\D$ is a couple $(u,\mu)$, where $u:\R_+\mapsto \R$
is a continuous function with $u(0)=0$ and whose distributional derivative
$u_x$ lies in $\L^2$, while $\mu=\mu^a+\mu^s$
is a bounded nonnegative Radon measure on $\R_+\,$,
whose absolutely continuous part (w.r.t.~Lebesgue measure) satisfies
\beq
d\mu^a = u_x^2\, dx\,.
\label{3.9}
\eeq
In the following, on the family of Radon measures on $\R_+$
we consider the distance
\beq 
d(\mu,\tilde\mu)\doteq \sup_\varphi \left|\int \varphi\,d\mu-
\int \varphi\,d\tilde \mu\right|\,,
\label{measdist}
\eeq
where the supremum is taken over all smooth functions $\varphi$
with $|\varphi|\leq 1$, $|\varphi_x|\leq 1$.

We recall that a semigroup $S$ on a domain $\D$ is a map
$S:\D\times [0,\infty[\,\mapsto\D$ such that
$S_0 w=w$ and $S_s(S_t w)=S_{s+t}w$ for every $s, t\geq 0$ and $w\in\D$.
\v

\begin{Theorem}{\label{Theorem3.1}}
Assume that the flux function $f$ satisfies condition (\ref{L0}).
Then there exists a semigroup $S:\D\times[0,\infty[\,\mapsto \D$
with the following properties.  Calling $t\mapsto 
S_t(\bar u,\bar \mu)=(u(t),\, \mu_{(t)})$
the trajectory corresponding to an initial data $(\bar u, \bar \mu)\in\D$,
one has:
\v
(i) The function $u=u(t,x)$ is  locally H\"older continuous in
$\R_+\times\R_+$. It provides a solution of (\ref{3.1})-(\ref{3.2})
in the sense of Definition 1.1 with initial condition
\beq
u(0,x)=\bar u(x)\,.
\label{3.10}
\eeq
\v
(ii) The assignment $t\mapsto \mu_{(t)}$ provides
a measure valued solution to the linear transport equation
\beq
w_t+\big[f'(u)w\big]_x=0\,,\qquad\qquad w(0)=\bar\mu\,.
\label{3.11}
\eeq
Moreover, the singular part of the measure
$f''(u(t))\cdot \mu_{(t)}$ vanishes at almost every time $t\geq 0$:
\beq
f''(u(t))\mu_{(t)}^s  = 0, \qquad \mbox{ a. e. } t.
\label{3.11a}
\eeq
\v
(iii) (Temporal continuity) 
 For every $M>0$, the above solution $u$ and the
corresponding measure $\mu$ satisfy the Lipschitz continuity property:
\beq
\int_0^M \left|u(t,x)-u(s, x)\right|\,\,dx \leq C|t-s|, 
\label{3.12}
\eeq
$$
d( \mu_{(t)}\,,~ \mu_{(s)} ) \leq C|t-s|,
$$
where the constant $C$ depends only on $M$, on the flux function $f$,
and on the total energy $\bar\mu(\R_+)< \infty$. 

\v
(iv) (Continuous dependence on the initial data) 
Finally, consider a sequence of initial conditions 
$(\bar u^n, \bar \mu^n) \in\D$ with $\bar u^n\to\bar u$ uniformly
on bounded sets and $d(\bar \mu^n,\bar \mu)\to 0$ as $n\to\infty$,
for some $(\bar u,\bar\mu)\in\D$. Then the corresponding solutions satisfy
\beq
u^n(t,x)\to u(t,x) \label{3.13}
\eeq
uniformly for $t,x$ in bounded sets, while
\beq
d(\mu^n_{(t)},\mu_{(t)})\to 0 \label{3.14}
\eeq
for every $t>0$.
\end{Theorem}

\v

\noindent{\bf Proof.} We treat here the case where $\bar \mu$ 
has compact support, say contained in the interval $[0,R]$,
so that $\bar u$ is constant for $x>R$.
The general case follows by an easy approximation argument.
The proof will be given in several steps. 
\v

\noindent{\bf 1. Construction of the solution.}
Let an initial data $(\bar u,\bar\mu)\in\D$ be given.
Set
$\bar \xi \doteq \bar \mu(\R_+) < \infty$.
On the semi-infinite strip 
$\left\{ t\geq 0\,,~\xi\in [0,\bar \xi]\,\right\}$ we
construct a function $U=U(t,\xi)$ by first setting
\beq
U(0, \xi) = \ov U (\xi)\doteq \bar u\left(\bar y(\xi)\right),
\label{3.15}
\eeq
where
\beq
\bar y(\xi) \doteq  
\inf \Big\{x\geq 0\ ;~~ {\bar \mu}\left([0, x]\right) \geq \xi \Big\}
\label{3.16}
\eeq
for $0<\xi\leq \bar \xi$,
while
\beq\label{y0}
\bar y(0)=\sup\{ x\ ;~~  {\bar \mu}\left([0, x]\right) =0\}\,.
\eeq
Observe that the map  $\xi\mapsto\bar y(\xi)$ is nondecreasing and 
left continuous, but it may well have upward jumps.
The provision (\ref{y0}) makes it continuous at the point $\xi=0$.  
In any case, the composed function
$\xi\mapsto \bar u ( \bar y ( \xi))$ is always continuous.  
For positive times,
the function $U$ is then defined to be the solution of
\beq
\f{\p U}{\p t}(t, \xi) = \f12 \int_0^\xi f''(U(t, \eta))\,d\eta 
\label{3.17}
\eeq
with initial data (\ref{3.15}). By the assumption of Lipschitz 
continuity of $f''$, the function $U$
can be obtained as the unique fixed point of a contractive
transformation. Details will be given at the next step.

Having constructed $U(t, \xi)$, the characteristic curves 
are obtained by solving the equation
\beq
y(0, \xi) = \bar y(\xi)\,,\qquad\qquad
\f{\p y}{\p t}(t,\xi)  = f'(U(t, \xi))\,.
\label{3.18}
\eeq
Explicitly:
\beq
y(t, \xi) = \bar y(\xi) + \int_0^tf'(U(\tau, \xi))\, d\tau. 
\label{3.19}
\eeq

Notice that $t\mapsto U(t,\xi)$ yields the values of our solution $u$
along the characteristic curve
$t\mapsto y(t,\xi)$ starting at $\bar y(\xi)$. A remarkable
feature of equation (\ref{3.1}) is that,
if the energy is conserved, then these values can be
determined in advance, before computing the actual position of the
characteristic curve. The image of the mapping
$$
\xi \to ( y(t, \xi), U(t, \xi))
$$
is now contained inside the graph of the desired solution $u(t,\cdot)$.
More precisely, for any given $(t,x)$ we define
\beq
u(t,x)=U(t,\xi(t,x)),\label{3.20}
\eeq
where
$$
\xi(t,x)\doteq \sup\left\{ \xi\,;~~y(t,\xi)\leq x\right\}\,.
$$
Finally, at time $t$  the corresponding measure $\mu_{(t)}$
is defined as the push-forward of the Lebesgue measure
on $[0,\bar\xi\,]\,$ through the mapping $\xi\mapsto y(t,\xi)$.
For each Borel set $J\subset \R_+$ we thus define
\beq
\mu_{(t)} (J) \doteq \hbox{meas}~\Big\{
\xi\in [0,\bar\xi\,]~; ~~~y(t,\xi)\in J\Big\}\,.
\label{3.21}
\eeq
\v

\noindent{\bf 2. A contractive transformation.}
Consider the space of continuous functions
$\C([0,\infty[\,\times[0,\bar\xi\,])$, with the equivalent weighted norm
\beq
\|U\|_*\doteq\sup_{t\geq 0,~\xi\in [0,\bar\xi\,]}
\,e^{-L\bar \xi t}\left|U(t,\xi)\right|\,,
\label{wn}
\eeq
where $L$ is a Lipschitz constant for the function $f''$.
The transformation $U\mapsto\T U$ is defined as
\beq
\T U\,(t,\xi)\doteq\bar u(\bar y(\xi))+ \f12 \int_0^t
\int_0^\xi f''(U(s, \eta))\,d\eta\,ds\, .
\label{3.23}
\eeq
If now $\|U-V\|_*=\delta$, then
$$
\Big|f''(U(\tau,\eta))- f''(V(\tau,\eta))\Big|\leq L
|U(\tau,\eta)-V(\tau,\eta)|\leq L \,e^{L\bar\xi \tau} \delta\,.
$$
For every $t\geq 0$ and $\xi\in [0,\bar\xi\,]$ we thus have
$$
\left|(\T U-\T V)(t,\xi)\right| \leq \f12
\int_0^t \left[\int_0^\xi L \,e^{L\bar\xi \tau}\,\delta \,d\eta\right]
\,d\tau \leq  \f12 \int_0^t L\xi \,e^{L\bar\xi \tau}\,\delta \,d\tau <
 \f12 e^{L\bar\xi t}\,\delta \,.
$$
By the above inequality we conclude
$$
\|\T U-\T V\|_*\leq \f12 \|U-V\|_*\,,
$$
proving the contractivity of the map $\T$.
By the contraction mapping theorem, it admits a unique fixed point
$U=U(t,\xi)$, defined on $\R_+\times [0,\bar\xi\,]\,$.
In turn, the function $u=u(t,x)$ and the measures $\mu_{(t)}$
are well defined by (\ref{3.20})-(\ref{3.21}).

\v
\noindent{\bf 3. Absolute continuity.}
We prove here that the map $\xi\mapsto U(t,\xi)$ is absolutely
continuous, for each $t\geq 0$.
Indeed, consider first the case $t=0$.
Let $[\xi_k,\xi'_k]$, with  $k=1,\ldots,N$, be  disjoint 
intervals contained in $[0,\bar\xi\,]$.
Assume that
$$
\sum_k |\xi_k'-\xi_k|\leq\ve\,.
$$
Call $I_1$ the set of indices $k$ such that
$$
\frac{|U(\xi_k')-U(\xi_k)|}{y(\xi'_k)-y(\xi_k)}\leq \sqrt\ve
$$
and let $I_2$ be the set of indices where
the above quantity is $>\sqrt\ve$.
Then
$$\sum_{k\in I_1} |U(\xi'_k)-U(\xi_k)|\leq 
\sqrt\ve\cdot \sum_{k\in I_1} | y(\xi'_k)-y(\xi_k)|
\leq \sqrt\ve\,R\,,$$
while
$$\sum_{k\in I_2} |U(\xi'_k)-U(\xi_k)|\leq \frac{1}
{\sqrt\ve}\sum_{k\in I_2} \frac{|U(\xi'_k)-U(\xi_k)|^2}
{y(\xi_k')-y(\xi_k)}
\leq  \frac{1}
{\sqrt\ve}\cdot \sum_{k\in I_2} \int_{y(\xi_k)}^{y(\xi'_k)}
u_x^2\,dx$$
$$\leq \frac{1}
{\sqrt\ve} \sum_{k\in I_2} |\xi'_k-\xi_k|\leq \sqrt\ve\,.$$
Together, the two above inequalities yield
$$\sum_{k=1}^N |U(\xi'_k)-U(\xi_k)|\leq (1+R)\,\sqrt\ve\,,$$
proving the absolute continuity of the map 
$\xi\mapsto U(0,\xi)$.

For $t>0$, the absolute continuity of $U(t,\cdot)$ follows from 
the absolute continuity of $U(0,\cdot)$ together with (\ref{3.17}). Indeed,
$$|U(t,\xi')-U(t,\xi)|\leq |U(0,\xi')-U(0,\xi)|+|\xi'-\xi|\cdot 
\frac{t}{2}\,\sup_u |f''(u)|\,.$$

As a consequence, the partial derivative
$U_\xi$ exists at a.e.~$(t,\xi)$. 
By (\ref{3.17}), it satisfies the evolution equation
\beq
\f{\p}{\p t} U_\xi (t,\xi)= \f12 f''(U(t,\xi))\,.
\label{2004a}
\eeq

On the other hand, the map $\xi\mapsto y(t,\xi)$ can be
discontinuous. However, if 
$$
\lim_{\xi\to \xi^*-} y(t,\xi)=y_1<y_2=\lim_{\xi\to \xi^*+} y(t,\xi)\,,
$$
then the function $u(t,\cdot)$ must be constant on the interval
$[y_1,y_2]$.

\v
\noindent{\bf 4. Measure transformations.}
To proceed, we first need to analyse
the regular and the singular part of the push-forward 
of Lebesgue measure, under a 
continuous non-decreasing transformation.
\v
\noindent {\bf Lemma 1.} {\it Let $U:[0,\bar\xi]\mapsto \R$ be 
absolutely continuous
with square integrable derivative.  Let $\xi\mapsto y(\xi)$ be such that
\beq
y(\xi)=y(0)+\int_0^\xi U_\xi^2(\zeta)\,d\zeta\,.
\label{yac}
\eeq
For $x\in [y(0), y(\bar\xi)]$ define the function $u=u(x)$ implicitly by 
\beq
 u(y(\xi))\doteq U(\xi)\,.
\label{yac1}
\eeq 
Let $\mu$ be the push-forward of Lebesgue measure
on $[0,\bar\xi]$ through the map $y$, i.e.
\beq 
\mu(J)\doteq \meas \{\xi\in [0,\bar \xi]\,;~~y(\xi)\in J\}\,.
\label{yac2}
\eeq

Then the absolutely continuous and the singular part of $\mu$ 
w.r.t.~Lebesgue measure are respectively given by
\beq\label{mua}
 \mu^a(A)=\meas \{\xi\in [0,\bar \xi]\,;~~y(\xi)\in A\,,~~U_\xi(\xi)
\not= 0\}\,.
\eeq
\beq 
\label{mus}
\mu^s(A)=\meas \{\xi\in [0,\bar \xi]\,;~~y(\xi)\in A\,,~~U_\xi(\xi)
= 0\}\,.
\eeq
In addition, on the set $[y(0),\, y(\bar\xi)]$ one has
\beq\label{mud}
d\mu^a= u_x^2\,dx\,.
\eeq
Viceversa, if both $U$ 
and the map $y$ are absolutely continuous and (\ref{yac1}), 
(\ref{yac2}), (\ref{mud}) are valid,
then (\ref{yac}) must hold.}

\v
\noindent {\bf Proof.} 
By (\ref{yac}), the image of a set $I\subseteq [0,\bar\xi]$ 
$$I_\ve\doteq \{\xi\in [0,\bar\xi\,]\,;~~|U_\xi(\xi)|\leq\ve\}$$
under the mapping $\xi\mapsto y(\xi)$ has Lebesgue measure
$$\meas(y(I))=\int_I U_\xi^2 (\xi)\,d\xi\,.$$
It is thus clear that the singular part of $\mu$
is the push-forward of Lebesgue measure on the set  
$$I_0\doteq \{\xi\in [0,\bar\xi\,]\,;~~U_\xi(\xi)=0\}$$
Next, for any fixed $\ve>0$ 
take a measurable set $J\subset[0,\bar\xi\,]$ such that
$$U^2_\xi(\xi)\geq\ve\qquad\qquad\hbox{for all}~\xi\in J\,.$$
Then 
$$\int_{y(J)} u_x^2(x)\,dx=\int_J\left[ U_\xi\,\frac{d\xi}{dy}
\right]^2\, \frac{dy}{d\xi}\cdot d\xi
=\int_J\left[ U_\xi\,U_\xi^{-2}
\right]^2\, U_\xi^2\cdot d\xi=\meas(J)\,.$$
Since $\ve>0$ was arbitrary, this proves (\ref{mua}), (\ref{mud}).
To prove the last statement, assume (\ref{yac1}), (\ref{yac2}) and (\ref{mud}).
Call 
$$J_\ve\doteq\{ \xi\in [0,\bar\xi\,]\,;~~
y_\xi(\xi)\geq\ve\}\,.$$
Observe that, for $\xi\in J_\ve$. the chain rule yields
\beq \label{chr}
u_x(y(\xi))\,y_\xi(\xi)= U_\xi(\xi)\,.
\eeq
For $0<a<b<\bar\xi$ we now obtain
\beq 
\int_{[y(a),y(b)]\cap y(J_\ve)} u_x^2(x)\,dx
= \int_{[a,b]\cap J_\ve} u_x^2(y(\xi))\,y_\xi(\xi)\,d\xi
= \meas( [a,b]\cap J_\ve)
\label{chr2}
\eeq
Since $a,b$ were arbitrary, this implies
\beq 
y_\xi(\xi)=[u_x^2(y(\xi))]^{-1}
\label{chr3}
\eeq
for $\xi\in J_\ve$. Together with (\ref{chr}) this yields
\beq\label{uder}
u_x(y(\xi))=U_{\xi}^{-1}(\xi)\,,\qquad \qquad 
y_\xi(\xi)=U_\xi^2(\xi)
\eeq
for all $\xi\in J_\ve$. Since $\ve>0$ is arbitrary, we conclude
$$
y(\xi) = y(0)+\int_0^\xi y_\xi(\zeta)\,d\zeta
       = y(0)+\lim_{\ve\to 0} \int_{[0,\xi]\cap J_\ve} y_\xi(\zeta)\,d\zeta
       = y(0)+\int_0^\xi U_\xi^2(\zeta)\,d\zeta\,,
$$
proving (\ref{yac}).
\v

\noindent{\bf 5. A class of regular solutions.}
Having constructed the trajectory $t\mapsto (u(t,\cdot),\,\mu_{(t)})$,
we still need to prove that it satisfies equation (\ref{3.1}), 
coupled with the initial and boundary conditions (\ref{3.10}), (\ref{3.2}).
 We carry out the analysis first assuming that the map
$\xi\mapsto \bar y(\xi)$ is absolutely continuous.
At the end, this assumption will be removed.

For each $t\geq 0$ and $\xi\in [0,\bar\xi\,]$ define
\beq\label{ydef2}
y(t, 0)  = \bar y(0)+t\,f'(0)\,, \qquad\qquad 
y(t,\xi)=y(t, 0)+\int_0^\xi U_\xi^2(t,\zeta)\,d\zeta\,.
\eeq
By (\ref{2004a}) this implies
\beq \label{yxit}
\f{\p}{\p t} y_\xi(t,\xi)=\f{\p}{\p t}
U_\xi^2(t,\xi)= f''(U(t,\xi))\,U_\xi(t,\xi)\,.
\eeq 
We now check that the function $y=y(t,\xi)$
defined at (\ref{ydef2}) coincides with the one defined at 
(\ref{3.19}). Indeed, by the second part of Lemma 1,
their derivatives $y_\xi$ coincide at time $t=0$ and 
satisfy the same equation (\ref{yxit}).
In particular, from (\ref{ydef2}) it is clear that the map $t\mapsto
y(t,\xi)$ is non-decreasing.  In particular, 
characteristics never cross each
other.

We begin by observing that the boundary condition (\ref{3.2})
is clearly satisfied, because
$$
u(t, 0)= U(t, 0)=U(0,0)+\int_0^t U_t(\tau,0)\,d\tau = 0\,.
$$
Moreover, the initial condition (\ref{3.10}) holds because of
the definitions (\ref{3.15})-(\ref{3.16}).

To check that the limit function $u$ satisfies (\ref{3.1}), 
fix a time $t>0$.   
Since $u(t,x)\equiv 0$ for $x\in [0,y(t, 0)]$, in this region
the equation (\ref{3.1}) trivially holds.
For almost every $x\in [y(t, 0), ~y(t,\bar\xi)]$, 
there exists a unique $\xi\in [0,\bar\xi\,]$ such that
$x=y(t,\xi)$. In this case, our construction yields
$$
u_t+f'(u)u_x = U_t(t,\xi) = \f12 \int_0^\xi f''(U(t,\zeta))\,d\zeta
= \f12 \int_0^{y(t,\xi)} f''(u(t,\cdot))\,d\mu_{(t)}\,.
$$
This implies (\ref{3.1}), provided that we can show the
identity of measures
\beq 
f''(u)\,u^2_x\,dx = f''(u)\,d\mu_{(t)}
\label{3.28a}
\eeq
for almost every time $t\geq 0$.
We shall now work toward a proof of (\ref{3.28a}).

Since the function $u$ is continuous, by covering the open  region
$$
\Big\{(t,x)\in \R_+\times\R_+\,;~~f''(u(t,x))\not= 0\Big\}
$$
with countably many sets of the form
$$
\Gamma\doteq \Big\{ (t,x)\,;~~t\in [t_1,t_2]\,,~~x\in [y(t,a),~ y(t, b)]
\Big\}\,
$$
it suffices to prove the following statement.

Assume that 
$$
f''(u(t,x))>\delta>0\qquad\qquad (t,x)\in \Gamma\,.
$$
Then, for a.e.~$t\in [t_1, t_2]$, 
the restriction of the measure $\mu_{(t)}$ to the interval
$[y(t,a),~y(t, b)]$ is absolutely continuous w.r.t.~Lebesgue measure
and satisfies $d\mu_{(t)}=u_x^2\,dx$.

By construction, as long as $U$ ranges in a region where $f''>\delta$
we have
\beq
\f{\p}{\p t} U_\xi(t,\xi) > \f{\delta}{2}\,.
\label{3.29}
\eeq
Hence, for any $\ve>0$,
$$
\meas\Big(\{ (t,\xi)\in\Gamma \,;~~|U_\xi(t,\xi)|<\ve\}\Big)<
\f{4\bar\xi}{\delta}\,\ve\,.
$$
Since $\ve>0$ here is arbitrary, we conclude that there exists
a set of times $\N$ of measure zero such that
$$
\meas\Big(\{ \xi\in [a,b]\,;~~U_\xi(t,\xi)=0\}\Big)=0
$$
for all times $t\notin \N$.
By Lemma 1, $t\notin\N$thus implies that the restriction of $\mu_{(t)}$
to the interval $[y(t,a), ~y(t, b)]$ is absolutely continuous 
w.r.t.~Lebesgue measure.
Furthermore, by (\ref{ydef2}), Lemma 1 
shows that  its density is $d\mu_{(t)}=u_x^2(t)\,dx$.
This concludes the proof of (i) and (ii) in Theorem 3.1,
at least in the case where the function $\xi\mapsto \bar y(\xi)$
is absolutely continuous.
\v
\noindent{\bf 6. General initial data.}
We now consider a more general initial data 
$(\bar u,\bar\mu)\in\D$.
In this case, the map
$\xi\mapsto \bar y(\xi)$ is non-decreasing, left continuous
but not necessarily continuous.
Its distributional derivative is thus a measure, say
$D_\xi\bar y=\sigma=\sigma^a+\sigma^s$.
By the assumptions, the absolutely continuous part satisfies
$$d\sigma^a=\overline U^2_\xi\,d\xi\,,$$
so that
$$\bar y(\xi)=\bar y(0)+\int_0^\xi \ov U^2_\xi (\zeta)\,d\zeta+
\sigma^s([0,\xi[\,)\,.$$
Consider a new initial condition $(\bar u^*,\bar \mu^*)$
defined by setting
$$\bar y^*(\xi)=\bar y(0)+\int_0^\xi \ov U^2_\xi (\zeta)\,d\zeta\,,
\qquad\qquad \bar u^*(\bar y^*(\xi))= \ov U(\xi)$$
$$\bar\mu^*(J)=\meas\{ \xi\,;~~\bar y^*(\xi)\in J\}\,.$$
By construction, for this new initial data the function
$\xi\mapsto y^*(0,\xi)=\bar y^*(\xi)$ is absolutely continuous.
Hence, by the previous analysis,
the corresponding function
$u^*(t,x)$ provides a solution to the initial-boundary value problem
(\ref{3.1})-(\ref{3.2}) with initial data $(\bar u^*,\bar \mu^*)$.
It is now easy to check that the function constructed in 
(\ref{3.17})--(\ref{3.20}) for the original initial data $\bar u$ satisfies
$$
u\Big(t,~y(t,\xi)+\sigma^s([0,\xi[)\Big)=U(t,\xi)\,.
$$
More precisely,
$$
u(t,x) = U(t,\xi) \quad\qquad\hbox{where} \quad
\xi = \inf\,\{\zeta\,; ~~y(t,\zeta)+ \sigma^s([0,\zeta])\geq x\}\,.
$$
By the previous analysis, $u^*$ provides a solution. 
Hence the same is true of $u$.

\v
\noindent{\bf 7. Continuity properties.}
Recall that $\bar \xi = \bar\mu(\R_+) < \infty$ is the total mass of each of 
the measures $\mu_{(t)}$. We have 
$$
\tv\{ u(t,\cdot)\,;~~[0, M]\,\}\leq \sqrt {\bar\xi M}\,.
$$
Since $u(t, 0)=0$, for any $x\in [0, M]$ we have
$$
|u(t,x)|\leq \sqrt {\bar\xi M}\, .
$$
This implies the Lipschitz continuity property w.r.t.~time:
\beq
\int_0^M |u(t,x)-u(s, x)|\,dx
\leq |t-s|\cdot \left\{\sup_{\omega}
|f'(\omega)|\cdot \sqrt { \bar\xi M}+\frac{\bar\xi M}{2}\cdot
\sup_{\omega} |f''(\omega)| \right\}
\label{3.32}
\eeq
where both sup are taken over $|\omega |\leq \sqrt {\bar\xi M}$.

Next, consider a convergent sequence of initial data 
$(\bar u^m,\bar\mu^m)_{m\geq 0}$.
The assumption of Theorem \ref{Theorem3.1} implies that the corresponding
functions $\bar U^n$ satisfy
$$
\ov U^m(\xi)\to \bar U(\xi)
$$
uniformly on $[0,\bar\xi]$. 
Therefore $U^m(t,\xi)\to U(t,\xi)$ uniformly on the
domain
$[0,T]\times [0,\bar\xi\,]\,$, for any $T>0$. 
In turn, this implies the convergence
(\ref{3.13})-(\ref{3.14}).

\v
\noindent{\bf 8. H\"older continuity.} We show that $u(t,x)$ is
H\"older continuous locally in $(t, x)$. First we know by Sobolev
embedding that $u$ is H\"older continuous in $x$ for each fixed time
$t$ with exponent $\alpha = 1/2$. In the time direction, we know that
the derivative of $u$ along a characteristic is bounded, thus $u$
is Lipschitz continuous in time along a characteristic. The characteristic
speed is $u$ which is locally bounded, thus the distance traveled in the
$x$ direction is order one of time. Combining the two parts, we conclude
that $u$ is H\"older continuous locally in both space and time.

This completes the proof of Theorem 3.1.
\v

\noindent{\bf Remark.} 
The previous construction of solutions to (\ref{3.1})-(\ref{3.2}) 
works equally well for negative times. The semigroup $S$
can thus be extended to a group $\Psi:\D\times\R\mapsto \D$.

\bigskip

\setcounter{equation}{0}
\section{Characterization of semigroup trajectories}

In the previous section, a solution $u$ to the initial-boundary value 
problem (\ref{3.1})-(\ref{3.2}), (\ref{3.10}),
was obtained as the fixed point of a contractive transformation.
Hence, any other solution which provides a fixed point to the same
transformation necessarily coincides with $u$.  
A straightforward uniqueness result can be stated as follows.

\begin{Theorem}{\label{Theorem4.1}}
Assume that $f$ satisfies (\ref{L0}). 
Consider a function $u=u(t,x)$ 
and a family of measures $\mu_{(t)}$ satisfying (i) and (ii)
in Theorem 3.1.  Moreover, calling
\beq 
y(t,\xi)\doteq 
\inf \Big\{x\geq 0\ ;~~ \mu_{(t)}\left([0, x]\right) \geq \xi \Big\}\,,
\label{4.1}
\eeq
\beq 
U(t, \xi) \doteq u(t,\, y(t,\xi)),
\label{4.2}
\eeq
assume that for a.e.~$\xi$ the map $t\mapsto U(t,\xi)$ is absolutely
continuous and satisfies the differential equation (\ref{3.17}).   Then
one has the identity
\beq 
(u(t), \mu_{(t)})= S_t(\bar u,\bar\mu).
\label{4.3}
\eeq
In particular, the solution which satisfies the above conditions 
is unique.
\end{Theorem}
\v

We conjecture that a uniqueness result remains valid
even without the assumption (\ref{3.17}) on the corresponding function $U$.
The basic ingredient toward a uniqueness result
is the assumption
\beq
 f''(u)\, d\mu^a{(t)}= f''(u) \,u_x^2(t)\, dx\,.
\label{4.4}
\eeq
for a.e.~$t$.
We now show that this is indeed the case
under the additional condition $f''>0$.

\begin{Theorem}{\label{Theorem4.2}}
In addition to assumption (\ref{L0}), let $f''(\cdot)>0$.
Consider a function $u=u(t,x)$ and a family of measures $\mu_{(t)}$ 
satisfying (i) and (ii) in Theorem 3.1. Then identity (\ref{4.3}) holds.
\end{Theorem}

Indeed, observe that the flow on $\L^1([0,\bar\xi])$ generated by the 
evolution equation (\ref{3.17}) 
is Lipschitz continuous w.r.t.~time and to the initial data.
Adopting a semigroup notation, call $t\mapsto V(t)=\S_t\ov V$
the trajectory corresponding to the initial data $\ov V\in 
\L^1([0,\bar\xi])$.
Since the couple $(u(t),\,\mu_{(t)})$ can be entirely recovered
from the function $U(t,\cdot)$ and the initial mapping
$\xi\mapsto \bar y(\xi)$, 
to prove uniqueness, it thus suffices to show that
\beq
\lim_{h\to 0+} \frac{1}{h} \int_0^{\bar\xi}
\Big|U(t+h,\xi)-(\S_h U(t))(\xi)\Big|\,d\xi=0 
\label{4.5}
\eeq
for almost every time $t>0$ (see Theorem 2.9 in \cite{[B1]}).
Since $f''>0$, our assumption implies that the singular part
of $\mu_{(t)}$ vanishes at a.e.~$t$.   
Choose a time $t$ where $\mu_{(t)}^s=0$.   Then
\beq 
U_\xi(t,\xi)\not= 0\qquad\hbox{ for a.e.}\quad\xi\in [0,\bar\xi\,]\,.
\label{4.6}
\eeq
Consider the map $\xi\mapsto y(t,\xi)$. Since 
$$
u_x^2\,dy = d\xi\,,\qquad\qquad 
U_\xi     = u_x\cdot \frac{dy}{d\xi} = \f1{u_x}\,,
$$
by (\ref{4.6}) the pre-image of a set of measure 
zero through the map $\xi\mapsto y(t,\xi)$ has measure zero.

If now $u=u(t,x)$ is differentiable at the point $(t, \,y(t,\xi))$,
we have the identity
\beq 
\begin{array}{rcl}
\f{\p}{\p t} U(t,\xi) & = & \Big[u_t +f'(u)u_x\Big](t,\, y(t,\xi))\\
                      & = & \f12\int_0^{y(t,\xi)} f''(u)\,u_x^2(t,x)\,dx
=\f12 \int_0^\xi f''(U(t,\eta))\,d\eta\,.
\end{array}
\label{4.7}
\eeq 

Observing that $u(t,\cdot)$ is differentiable at a.e.~$x$,
we conclude that (\ref{4.7}) holds at a.e.~$\xi\in[0,\bar\xi\,]$.
In turn, this implies (\ref{4.5}),
proving the theorem.

\endproof

\medskip

Notice how the condition on the vanishing of the singular part
is essential to ensure uniqueness. Otherwise, in Example 1
the solution $u(t,x)\equiv 0$ for $t\geq 1$, with
$\mu_{(t)}$ containing a unit mass at the origin, would
satisfy all the other requirements of the theorem.

\bigskip

\setcounter{equation}{0}
\section{A semigroup of dissipative solutions}

Next, we  examine dissipative solutions.
A major difference with the conservative case is that here the Cauchy
problem is well-posed if the flux function $f$ is strictly convex, but
ill posed otherwise, as shown in the next section.

In this section, our main concern will be the construction of a semigroup
of dissipative solutions under the additional assumption that
$f''\geq 0$. As domain $\D$ of our semigroup we choose
the space
$$\D\doteq \Big\{u:\R_+\mapsto \R\,,~~~u~\hbox{is 
absolutely continuous,}~~
u(0)=0\,,~~u_x\in \L^2\Big\}\,.$$

\begin{Theorem}{\label{Theorem5.1}}
Assume that the flux function $f$ satisfies (\ref{L0}) and $f'' \geq 0$.
Then there exists a semigroup $S:\D\times[0,\infty[\,\mapsto \D$
with the following properties.  Calling $t\mapsto u(t)= 
S_t\bar u$
the trajectory corresponding to an initial data $\bar u\in\D$,
one has:
\v
(i) The function $u=u(t,x)$ is  H\"older continuous. It
provides a solution of (\ref{3.1})-(\ref{3.2}) with
initial condition $u(0,x)=\bar u(x)$.
\v
(ii) For every $M>0$, the above solution $u$ 
satisfies the Lipschitz continuity property in time:
\beq
\int_0^M \left|u(t,x)-u(s, x)\right|\,\,dx \leq C|t-s|, 
\label{5.1}
\eeq
\v
(iii) Given a sequence of initial conditions 
$\bar u^n \in\D$, assume that 
$$
\|\bar u^n_x- \bar u_x\|_{\L^2([0,M])}\to 0
$$
for every $M>0$. 
Then the corresponding solutions satisfy
\beq
u^n(t,x)\to u(t,x) \label{5.2}
\eeq
uniformly for $t,x$ in bounded sets.
\end{Theorem}

\noindent {\bf Proof.}
Consider an initial condition $\bar u\in \D$.  For simplicity,
we again assume that $\bar u$ is constant outside a bounded interval, 
say $[0,R]$.
The general case follows from
an approximation argument.

To construct the
corresponding trajectory we begin by setting
$$\bar \xi\doteq \int_0^R|u_x^2(x)|\,dx\,.$$
Then we define the initial data
$$
\ov U(\xi)\doteq \bar u(\bar y(\xi))\, ,
$$
where
\beq 
\bar y(\xi)\doteq \inf\Big\{ x\geq 0\,;~~\int_0^x u_x^2(x)\,dx\geq \xi\Big\}\,.
\label{5.3}
\eeq
By the analysis in Section 3, the map $\xi \mapsto \ov U(\xi)$
is absolutely continuous, hence its derivative
$$\ov Z(\xi)=\frac{\partial}{\partial \xi}\ov U(\xi)$$
is a well defined function in $\L^1([0,\bar\xi\,])$.

Define the subset
$$
J^-\doteq \Big\{ \xi\in [0,\bar\xi\,]\,; ~~\ov Z(\xi)\leq 0\,\Big\}\,.
$$
Let $L$ be a Lipschitz constant for $f''$.
On the space of continuous functions $Y:\R_+\mapsto L^1([0,\bar\xi\,])$
with weighted norm
$$
\|Y\|_*\doteq \sup_t ~ e^{-L\bar \xi t} \|Y(t)\|_{L^1}\,,
$$
we now define a continuous transformation $Y\mapsto \T Y$ as follows.
\beq
\T Y(t,\xi)\doteq \ov Z(\xi)+\int_0^t \f12 f''\left(
\int_0^\xi \Phi(\eta, Y(s,\eta))\, d\eta\right) ds\, ,
\label{5.4}
\eeq
where
$$
\Phi(\eta,Y)=\min\{Y,0\}\qquad \hbox{if}~~\eta\in J^-\,,
$$
$$
\Phi(\eta,Y)=Y \qquad \hbox{if}~~\eta\in [0,\bar\xi\,]\setminus J^-\,.
$$
To check that $\T$ is a strict contraction,
assume that $\|Y-\Tilde Y\|_*=\kappa$, so that
$$
\int_0^{\bar \xi} |Y(t,\xi)-\Tilde Y(t,\xi)|\,
d\xi\leq \kappa\, e^{L\bar\xi t}
$$
for all $t\geq 0$. Then for every $s\geq 0$
$$
\int_0^\xi\Big|\Phi(\eta ,\,Y(s,\eta))- \Phi(\eta ,\,
\Tilde Y(s,\eta))\Big|\,d\eta\leq \kappa \,  e^{L\bar\xi s},
$$
and therefore
$$
\int_0^{\bar \xi} |(\T Y-\T \widetilde Y)(t,\xi)|d\xi\leq
\int_0^{\bar\xi}\int_0^t \frac{L\kappa}{2}  \, e^{L\bar\xi s}\, ds \,d\xi
\leq \frac{\kappa}{2}\, e^{L\bar\xi t}.
$$
By the definition of our weighted norm, this implies
$$\|\T Y-\T \widetilde Y\|_*\leq \f12 \,\|Y-\Tilde Y\|_*\,.$$
Let now $Y=Y(t,\xi)$ be the unique fixed point of the transformation $\T$.
Then one easily checks that the function
$$Z(t,\xi)\doteq Y(t,\xi)\qquad\hbox{if}~~\xi\notin J^-\,,$$
$$Z(t,\xi)\doteq \min\{ Y(t,\xi), ~~0\}\qquad\hbox{if}~~\xi\in J^-\,,$$
provides a solution to the equations
$$
Z(0,\xi)=\ov Z(\xi)\,,
$$
$$
\f{\p Z}{\p t}(t,\xi)=\f12 f''\left(\int_0^\xi
Z(t,\eta)\,d\eta\right)\qquad\hbox{if} ~~Z(s,\xi)\not= 0~~\hbox{for all}
~s\in [0,t]\,,
$$
$$ 
\f{\p Z}{\p t}(t,\xi)=0 \qquad\hbox{if} ~~Z(s,\xi)= 0~~\hbox{for some}
~s\in [0,t]\,.
$$
In turn, we can now define
$$
U(t,\xi)\doteq \int_0^\xi Z(t,\eta)\,d\eta
$$
and the characteristic curves
$$y(t,\xi)\doteq \bar y(\xi)+\int_0^t f'(U(s,\xi))\,ds\,.$$
In a similar way as in Section 3,
the dissipative solution  $u$ can now be obtained by setting
$$
u(t,x)=U(t,\xi(t,x)),
$$
where
$$\xi(t,x)\doteq \sup\left\{ \xi\,;~~y(t,\xi)\leq x\right\}\,.
$$

To see why  this construction actually yields a solution to 
(\ref{3.1}),
consider first the case where the map
$\xi\mapsto \bar y(\xi)$ is absolutely continuous.
Then
$y_\xi(0,\xi)=\ov U^2_\xi(0,\xi)=Z^2(0,\xi)$.
Since
$$
\frac{\partial}{\partial t} y_\xi=f''(U)\,U_\xi=f''(U)\,Z=
\frac{\partial}{\partial t} Z^2
$$
for all $t,\xi$ we deduce the identity
$$y_\xi(t,\xi)=Z^2(t,\xi)=U^2_\xi(t,\xi).$$
Moreover, (3.34) again holds.
As in the proof of Theorem 3.1, we obtain the relations
\beq\label{cvar}
Z(t,\xi)=\frac{1}{u_x(t,y(t,\xi))}\,,\qquad\qquad 
[y_\xi(t,\xi)]^{-1}=u_x^2(t,y(t,\xi)).
\eeq
For almost every $x\in [ y(t, 0)\,,~y(t,\bar\xi)]\,$, if 
$x=y(t,\xi)$, then 
$$[u_t+f'(u)u_x](t,x)=\frac{d}{dt} u(t, y(t,\xi))=
\frac{\partial}{\partial t} U(t,\xi)
=\int_0^\xi\frac{\partial}{\partial t} Z(t,\eta)\,d\eta$$
$$= \f12\int_0^\xi f''(U(t,\eta))\,d\eta
=\f12\int_0^x f''(u(t,y))\,u_x^2(t,y)\,dy\,.$$
The second identity in (\ref{cvar}) was used here to
change the variable of integration.

The extension to
the case of general initial data, where the map $\xi\mapsto y(t,\xi)$
is not necessarily absolutely continuous, is carried out as in
the earlier proof of Theorem 3.1.  We skip the details.

\bigskip

\setcounter{equation}{0}
\section{Instability of dissipative solutions for non-convex flux}
 
In this section, we show that if the convexity assumption $f''\geq 0$
is dropped, then the Cauchy problem for the equation (\ref{3.1})-(\ref{3.2})
is ill posed, in general.

\bigskip
\noindent{\bf Example 2.} Consider the flux function
$f(u)=u^3$.    Let $U=U(t,\xi)$ be a solution of
(\ref{3.17}), with $\xi\in [0,3]$, such that at some time $t_0>0$ there holds
$$
U(t_0,\xi) =
\left\{
\begin{array}{ll}
   \xi ,   & \qquad  \xi \in [0,1]\,, \nonumber \\
 2-\xi ,   & \qquad  \xi \in [1,2]\,, \nonumber\\
    0  ,   & \qquad  \xi \in [2,3]\,. \nonumber\\
\end{array}
\right.
$$
Consider first the conservative solution $u=u(t,x)$. This is well defined
forward and backward in time.
At time $t=t_0$, its explicit values are
$$
u(t_0,x)=\left\{
\begin{array}{ll}
  x\,,  & \qquad  x \in [0,1]\,, \nonumber \\
2-x\,,  & \qquad  x \in [1,2]\,, \nonumber\\
  0\,,  & \qquad  x >2 \, \nonumber
\end{array}
\right.
$$
while a unit mass is concentrated at the point $x=2$. 
Assuming $t_0$ sufficiently small,
we have
$$
U_t(t,\xi)=\int_0^\xi 3U(t,\eta)\,d\eta >0
$$
for all $t\in [0,t_0]$ and $\xi\in\,]0,3]$.
Hence
$$\frac{\partial}{\partial t} U_\xi(t,\xi)=3U <0,\qquad U_\xi(t,\xi)<0
\qquad\quad\hbox{for}~~~ t\in [0,t_0[\,,~~2<\xi<3\,. $$

Next, consider a dissipative solution $v$ coinciding with $u$ at time $t=0$.
This means
\beq 
v(0,x)=u(0,x)=U(0,\xi)\qquad\qquad \hbox{for}~~x=y(0,\xi)\,.\label{6.1}
\eeq
We recall that
$$y(t,\xi)
=\int_0^\xi U_\xi^2(t,\eta)\,d\eta\,.$$
Clearly, $v$ will still coincide with $u$ 
as long as its gradient remains bounded (equivalently,
as long as $U_\xi$ remains bounded away from zero).
On the other hand, for $t>t_0$,  
the dissipative solution $v=v(t,x)$ coincides with
the conservative one only on the interval where $x\leq
y(t, 2)$, while $v$ is constant for $x\geq y(t, 2)$. In other words,
$$v(t,x)= u(t,x)\qquad\qquad\hbox{if}~~t\in [0,t_0]\,,$$
$$v(t,x)=\left\{\begin{array}{ll}
u(t,x) & \qquad 0\leq x\leq y(t, 2)\,,\\
u(t, y(t, 2)) & \qquad x>y(t, 2)\,,
\end{array}\right.\qquad\hbox{if}~~t\in [t_0,\,2t_0]\,.
$$
Energy dissipation occurs at time $t=t_0$,
namely
$$
\int_0^\infty v_x^2(t,x)\,dx =  \left\{
\begin{array}{ll}
3 & \qquad t\in [0,t_0[,\\
2 & \qquad t \geq  t_0.
\end{array}
\right.
$$
Next, consider a family of perturbed initial conditions, say
$$U^\ve(0,\xi)=U(0,\xi)+\ve\phi(\xi)\,,$$
where $\phi$ is a non-negative smooth function, whose support is
contained in $[0, 1]$. 
Since $U\mapsto f''(U)=6U$ is
a monotone increasing function, by a comparison argument
from (\ref{3.17}) we deduce
$$U^\ve(t,\xi)\geq U(t,\xi)$$
for all $\ve,t>0$, $\xi\in [0,3]$.
In fact, for a nontrivial $\phi$ we can assume a strict inequality:
$$U^\ve(t,\xi)> U(t,\xi)\qquad\qquad t>0\, ,~~~\xi\in [2,3]\,.$$
For $2<\xi<3$ we now use the relations
$$
\frac{\partial}{\partial t} U^\ve_\xi(t,\xi)= 3U^\ve(t,\xi)
> 3U(t,\xi)=
\frac{\partial}{\partial t} U_\xi(t,\xi)\,,
\qquad\qquad U^\ve_\xi(0,\xi)=U_\xi(0,\xi)\,,
$$
and deduce
$$U_\xi^\ve (t,\xi)>U_\xi(t,\xi)\geq 0\qquad\qquad 
t\in [0,t_0]\,.$$
Moreover, for $t\geq t_0$ and $2<\xi<3$ one has
$$\frac{\partial}{\partial t} U^\ve_\xi(t,\xi)= 3U^\ve(t,\xi)>
3U(t,\xi)\geq 0\,.$$
Therefore, for each $\ve>0$, the quantity
$U^\ve_\xi(t,\xi)$ is still strictly positive 
at time $t=t_0$ and increases afterwards.  
It thus remains uniformly bounded away from zero.

Since $u_x=U_\xi^{-1}$,
the above implies that,
for any fixed $\ve>0$, the corresponding conservative
solution $u^\ve=u^\ve(t,x)$
has a uniformly bounded gradient.   The dissipative solution
thus coincides with the conservative one.
As $\ve\to 0$, at time $t=0$ our construction yields
$$\|u^\ve(0)-u(0)\|_{\C^0}\to 0\,,\qquad\qquad
\|u_x^\ve(0)-u_x(0)\|_{\L^2}\to 0\,.
$$
However, when $t>t_0$ and $x> y(t, 2)$ the previous analysis yields
$$\lim_{\ve\to 0+} u^\ve(t,x)= u(t,x)\not= v(t,x)\,,$$
where $u,v$ are respectively the conservative and the dissipative solutions
of (\ref{3.1})-(\ref{3.2}), with the same initial data (\ref{6.1}).
The example proves that
dissipative solutions do not depend continuously on the inital data.
\v

\noindent {\bf Remark.}
The previous example also shows that
the family of dissipative solutions
may not be closed.  Since the set of solutions which are limits of
vanishing viscosity approximations is closed and connected,
we see that this set cannot coincide with the set of dissipative solutions.
\bigskip

{\bf Acknowledgments.}  We thank John Hunter for helpful conversations.
This work is supported in part by NSF DMS-0305497 and NSF DMS-0305114
for Yuxi Zheng, NSF of China under Grants 10131050 and 10276036
(the innovation grants from Chinese Academy of Sciences) for Ping Zhang,
and by the Italian M.I.U.R. within the research project
\# 2002017219 ``Equazioni iperboliche e paraboliche non lineari'' 
for Alberto Bressan. 
This work was started when Ping Zhang visited Penn State 
University. He would like to thank the department for its hospitality.


\begin{thebibliography}{99}


\bibitem{[B1]} A. Bressan, 
{\it Hyperbolic Systems of Conservation Laws. 
The One-Dimensional Cauchy Problem.}
Oxford Univ. Press, 2000.

\bibitem{[BC]} A. Bressan and A.~Constantin, work in progress.

\bibitem{[BZ]} A. Bressan and Yuxi Zheng, 
Conservative solutions to a nonlinear variational wave equation, work in
progress.

\bibitem{[3]} R.~Camassa and D.~Holm, An integrable shallow water equation
with peaked solitons, {\it Phys. Rev. Lett.}, 71(1993), 1661--1664.

\bibitem{[4]} A.~Constantin, On the Cauchy problem for the periodic 
Camassa-Holm equation, {\it J. Diff. Equations}, {\bf 141}(1997), 218--235.

\bibitem{[5]} A.~Constantin and H.~P.~McKean, A shallow water equation on the
circle, {\it Comm. Pure Appl. Math.}, {\bf 52}(1999), 949--982.

\bibitem{[6]} A.~Constantin and J.~Escher, Wave breaking for nonlinear
nonlocal shallow water equations, {\it Acta Mathematica}, 
{\bf 181}(1998), 229--243.

\bibitem{[7]} A.~Constantin and J.~Escher, Global weak solutions for a
shallow water equation, {\it Indiana Univ. Math. J.},
{\bf  47}(1998),  1527--1545. 
 

\bibitem{[8]} A.~Constantin and L.~Molinet, Orbital stability of solitary waves
for a shallow water equation, {\it  Phys. D}, {\bf  157}(2001), 75--89. 


\bibitem{[13]} R.~T.~Glassey, J.~K.~Hunter and Yuxi Zheng, Singularities in a
Nonlinear Variational Wave equation, {\it J. Differential Equations},
{\bf 129}(1996), 49-78.

\bibitem{[14]} R.~T.~Glassey, J.~K.~Hunter and Yuxi Zheng, Singularities and
oscillations in a nonlinear variational wave equation, in
{\it Singularities and Oscillations}, edited by J.~Rauch and M.~Taylor,
IMA, {\bf 91},  Springer, 1997.

\bibitem{[15]} J.~K.~Hunter, and R.~A.~Saxton, Dynamics of director fields,
{\it SIAM J. Appl. Ma- th.}, {\bf 51}(1991), 1498-1521.

\bibitem{[16]} J.~K.~Hunter, and Yuxi Zheng, On a nonlinear hyperbolic
variational equation I and II, {\it Arch. Rat. Mech. Anal.},
{\bf 129}(1995), 305-353 and 355-383.


\bibitem{[20]} R.~A.~Saxton, Dynamic instability of the liquid crystal
director, in {\it Contemporary Mathematics}, Vol.~{\bf 100}:
Current Progress in Hyperbolic Systems,
pp.~325--330, ed. W.~B.~Lindquist, AMS, Providence, 1989.



\bibitem{[22]} Zhouping Xin and Ping Zhang, On the weak solutions to a shallow
water equation, {\it Comm. Pure Appl. Math.}  {\bf 53}(2000),  1411--1433.


\bibitem{[23]} Ping Zhang and Yuxi Zheng, On oscillations of an
asymptotic equation of a nonlinear variational wave equation,
{\it Asymptotic Analysis}, {\bf 18}(1998), 307--327 .

\bibitem{[24]} Ping Zhang and Yuxi Zheng, On the existence and
 uniqueness to an asymptotic equation of a variational wave
equation,  {\it Acta Mathematica Sinica}, {\bf 15}(1999), 115--130.

\bibitem{[25]} Ping Zhang and Yuxi Zheng, On the existence and
uniqueness to an asymptotic equation of a variational wave
equation with general data, {\it Arch. Rat. Mech. Anal.} 
{\bf 155}(2000), pp. 49--83.


\bibitem{[26]} Ping Zhang and Yuxi Zheng, Rarefactive solutions to a
nonlinear variational wave equation,  {\it Comm. Partial Differential
Equations}, {\bf 26}(2001), pp. 381-419.

\bibitem{[27]} Ping Zhang and Yuxi Zheng, Singular and rarefactive solutions
to a nonlinear variational wave equation, {\it Chinese Annals of
Mathematics}, {\bf 22B}, 2(2001), pp. 159-170.

\bibitem{[28]} Ping Zhang and Yuxi Zheng, Weak solutions to a nonlinear
variational wave equation, {\it Arch. Rat. Mech. Anal.}, 
{\bf 166} (2003), 303--319.


\bibitem{zzpV} Ping Zhang and Yuxi Zheng,  On the Second-Order 
Asymptotic Equation of a Variational Wave Equation, 
{\it Proc A of The Royal Soc Edinburgh, A. Mathematics}, 
{\bf 132A}(2002), 483--509. 

\bibitem{[29]} Ping Zhang and Yuxi Zheng, Weak Solutions to A
Nonlinear Variational Wave Equation with General Data,
Annals of Inst H. Poincar\'e (in press), 2004.

\end{thebibliography}
\end{document}